\documentclass[12pt]{article}

\usepackage{a4, amssymb}
\usepackage{a4wide}
\usepackage{amsmath}
\usepackage{amsfonts}
\usepackage{amssymb}
\usepackage{graphics}
\usepackage{textcomp}
\usepackage{latexsym}
\usepackage{epsfig}
\newtheorem{theorem}{Theorem}[section]
\newtheorem{proposition}[theorem]{Proposition}
\newtheorem{corollary}[theorem]{Corollary}
\newtheorem{lemma}[theorem]{Lemma}
\newcommand{\IR}{\ensuremath{\mathbb{R}}}

\newcommand{\IZ}{\ensuremath{\mathbb{Z}}}
\newcommand{\IN}{\ensuremath{\mathbb{N}}}
\newcommand{\IP}{\ensuremath{\mathbb{P}}}
\newcommand{\IE}{\ensuremath{\mathbb{E}}}
\newcommand{\Var}{{\mathbb{V}}{\rm ar}}
\def\QED{{\nobreak\hfil\penalty50\hskip1em\hbox{}
\nobreak\hfil$\Box$\parfillskip0pt\par}\bigskip}
\newcommand{\plg}{\ensuremath{p(\lambda,g)}}
\newcommand{\plgr}{\ensuremath{p(\lambda,g_R)}}
\newcommand{\plgrr}{\ensuremath{p(\lambda,g^R)}}
\newcommand{\plnndg}{\ensuremath{p(\lambda_n/n^d,g)}}
\newcommand{\plnndgr}{\ensuremath{p(\lambda_n/n^d,g_R)}}
\newcommand{\plnndgrr}{\ensuremath{p(\lambda_n/n^d,g^R)}}
\newcommand{\plngn}{\ensuremath{p(\lambda_n,g_n)}}
\newcommand{\plngrn}{\ensuremath{p(\lambda_n,g_{R,n})}}
\newcommand{\plngrrn}{\ensuremath{p(\lambda_n,g^{R,n})}}
\newcommand{\plgg}{\ensuremath{p_{\lambda}^{g,g}}}
\newcommand{\plgrg}{\ensuremath{p_{\lambda}^{g_R,g}}}
\newcommand{\plgrgr}{\ensuremath{p_{\lambda}^{g_R,g_R}}}
\newcommand{\plnndgg}{\ensuremath{p_{\lambda_n/n^d}^{g,g}}}
\newcommand{\plnndgrg}{\ensuremath{p_{\lambda_n/n^d}^{g_R,g}}}
\newcommand{\plnndgrgr}{\ensuremath{p_{\lambda_n/n^d}^{g_R,g_R}}}
\newcommand{\plngngn}{\ensuremath{p_{\lambda_n}^{g_n,g_n}}}
\newcommand{\plngrngn}{\ensuremath{p_{\lambda_n}^{g_{R,n},g_n}}}
\newcommand{\plngrngrn}{\ensuremath{p_{\lambda_n}^{g_{R,n},g_{R,n}}}}

\newcommand{\plngnr}{\ensuremath{p(\lambda_n,g_{n,R})}}
\newcommand{\plngnrr}{\ensuremath{p(\lambda_n,g_n^R)}}
\newcommand{\plnndhnr}{\ensuremath{p(\lambda_n/n^d,k_{nR})}}
\newcommand{\plnndhnnr}{\ensuremath{p(\lambda_n/n^d,k^{nR})}}
\newcommand{\plngnrgnr}{\ensuremath{p_{\lambda_n}^{g_{n,R},g_{n,R}}}}
\newcommand{\plngnrgn}{\ensuremath{p_{\lambda_n}^{g_{n,R},g_n}}}
\newcommand{\plnndhnrhnr}{\ensuremath{p_{\lambda_n/n^d}^{k_{nR},k_{nR}}}}
\newcommand{\plnndhnrg}{\ensuremath{p_{\lambda_n/n^d}^{k_{nR},g}}}
\newcommand{\Cov}{{\mathbb{C}}{\rm ov}}

\newcommand{\QEDnoskip}{\ensuremath{\hfill \Box}}
\newcommand{\AAA}{\ensuremath{\mathcal{A}}}

\newcommand{\FFF}{\ensuremath{\mathcal{F}}}

\newcommand{\grotebox}{\ensuremath{\begin{picture}(12,10) 
\put(1,-2){\framebox(10,10){}} \end{picture}}}
\newcommand{\verschovenbox}{\ensuremath{\begin{picture}(12,10) 
\put(1,-2){\line(1,0){6}} \put(1,-2){\line(0,1){6}} \put(5,2){\framebox(6,6){}} 
\end{picture}}}
\newcommand{\gebiedL}{\ensuremath{\begin{picture}(12,10) 
\put(1,-2){\line(1,0){10}} \put(1,-2){\line(0,1){10}} \put(1,8){\line(1,-1){10}} 
\end{picture}}}
\newcommand{\gebiedgeroteerdeL}{\ensuremath{\begin{picture}(12,10) 
\put(11,8){\line(-1,0){10}} \put(11,8){\line(0,-1){10}} 
\put(1,8){\line(1,-1){10}} \end{picture}}}

\begin{document}

\begin{titlepage}

\title{On central limit theorems in the random connection model}

\author{Tim van de Brug and Ronald Meester} 

\date{August 16, 2003}

\maketitle
\begin{center}
Divisie Wiskunde,
Faculteit der Exacte Wetenschappen,
Vrije Universiteit,
De Boelelaan 1081a,
1081 HV Amsterdam, The Netherlands, {\tt rmeester@cs.vu.nl}. 
\end{center}

\begin{abstract}
Consider a sequence of Poisson random 
connection 
models $(X_n,\lambda_n,g_n)$ on $\IR^d$, where $\lambda_n/n^d \to \lambda>0$ and 
$g_n(x)=g(nx)$ for some
non-increasing, integrable connection function $g$. Let $I_n(g)$ be the number 
of 
isolated vertices of 
$(X_n,\lambda_n,g_n)$ in some bounded Borel set $K$, where $K$ has non-empty 
interior 
and 
boundary of Lebesgue measure zero. Roy and Sarkar (2003) claim that
$$
\frac{I_n(g)-\IE I_n(g)}{\sqrt{\Var I_n(g)}} \rightsquigarrow N(0,1),\qquad 
n\to\infty,
$$
where $\rightsquigarrow$ denotes convergence in distribution. However, their 
proof has errors. We correct their proof and extend the result to larger
components when the connection function $g$ has bounded support.
\end{abstract}

\end{titlepage}

\section{Introduction}

Let $(X,\lambda,g)$ denote a Poisson random connection model, where $X$ is the 
underlying Poisson point process on $\IR^d$ with density $\lambda>0$, and where 
$g$ is a {\em connection function} which we assume to be a non-increasing and 
which satisfies 
$0 < \int_{\IR^d} g(|x|)\,dx < \infty$. In words, this amounts to saying that 
any two points $x$ and $y$ of $X$ are connected with probability $g(|x-y|)$, 
independently of all other pairs, independently of $X$. The random connection 
model plays an important role in many areas, for instance in telecommunications 
and epidemiology. In telecommunications, the points of the point process can 
represent base stations, and the connection function then tells us that two base 
stations at locations $x$ and $y$ respectively, can communicate to each other 
with probability $g(|x-y|)$. In epidemiology, the connection function can for 
instance represent the probability that an infected herd at location $x$ infects 
another herd at location $y$.

Let $K$ be a bounded 
Borel subset of 
$\IR^d$ with non-empty interior and boundary of Lebesgue measure zero. 
Consider a sequence of positive real numbers $\lambda_n$ with 
$\lambda_n/n^d\to\lambda$, let $X_n$ be a Poisson 
process on $\IR^d$ with density $\lambda_n$ and let $g_n$ be 
the connection function defined by $g_n(x)=g(nx)$. Consider the sequence of 
Poisson random 
connection models $(X_n,\lambda_n,g_n)$ on $\IR^d$. Let $I_n(g)$ be the number 
of 
isolated 
vertices of $(X_n,\lambda_n,g_n)$ in $K$. Roy and Sarkar (2003) claim to prove 
the 
following result.
\begin{theorem} \sl
\label{mh}
\begin{equation}\label{eq0}
\frac{I_n(g)-\IE I_n(g)}{\sqrt{\Var I_n(g)}} \rightsquigarrow N(0,1),\qquad 
n\to\infty,
\end{equation}
where\/\ $\rightsquigarrow$ denotes convergence in distribution. 
\end{theorem}
Although the statement of this result is correct, the
proof in Roy and Sarkar (2003) is not. In this note, we explain what went wrong 
in their proof, and how this
can be corrected. In addition, we prove an extension to larger components in 
case the connection function has bounded support.

\section{Truncation and scaling}

The central limit theorem (\ref{eq0}) is relatively easy to show when $g$ has
bounded support, see Roy and Sarkar (2003). Hence, the strategy adopted by Roy 
and 
Sarkar (2003) is to truncate the relevant connection functions, and let the 
truncation 
go 
to infinity. This means that there are two operations involved: scaling and 
truncation. The root of the problem lies in the fact that these two operations 
do not commute. 

Following Roy and Sarkar (2003), we define for $R>0$ and $n\in\IN$ connection 
functions $g_R,g^R,g_{n,R},g_n^R : 
[0,\infty)\to [0,1]$ by
\[
g_R(x)=1_{\{x\leq R\}}g(x),~ g^R(x)=1_{\{x>R\}}g(x),~ g_{n,R}(x)=1_{\{x\leq 
R\}}g(nx),~ g_n^R(x)=1_{\{x>R\}}g(nx),
\]
where the {\em indicator function} $1_{x \leq R}$ is by definition equal to 1 
when $x \leq R$ and equal to 0 when $x > R$, and similarly for the other 
indicator functions.
Note that the notation $g_R$ can formally not be used to 
denote $1_{\{\cdot\leq R\}} g(\cdot)$, since $g_n$ has already been defined as 
$g(n\cdot)$. Nevertheless we shall adopt this notation, because we think that 
this 
will not cause any confusion. Henceforth $g_R$ will always denote 
$1_{\{\cdot\leq 
R\}} g(\cdot)$ and $g_n$ will always denote $g(n\cdot)$. Let $L_R(g)$ be the 
number 
of isolated vertices of $(X,\lambda,g_R)$ in $K$ that are not isolated in 
$(X,\lambda,g)$. Let $J_{n,R}(g)$ be the number of isolated vertices of 
$(X_n,\lambda_n,g_{n,R})$ in $K$ and let $L_{n,R}(g)=J_{n,R}(g)-I_n(g)$ be the 
number of isolated vertices of $(X_n,\lambda_n,g_{n,R})$ in $K$ that are not 
isolated in $(X_n,\lambda_n,g_n)$. 

Roy and Sarkar (2003) claim the following (without proof).

\medskip\noindent
{\bf Statement A}\, {\sl If (\ref{eq0}) is true when the connection function\/\ 
$g$
has bounded support, then it is the case that 
\begin{equation}\label{eq0l}
\frac{J_{n,R}(g)-\IE J_{n,R}(g)}{\sqrt{\Var J_{n,R}(g)}} \rightsquigarrow 
N(0,1),\qquad n\to\infty,
\end{equation}
for\/\ {\em any} connection function\/\ $g$.}

\medskip\noindent
They then proceed, via a number of moment estimates involving $J_{n,R}(g)$ and 
$L_{n,R}(g)$, to show that the truth of (\ref{eq0l}) for any connection function 
$g$, implies the full central limit theorem in (\ref{eq0}). 

One problem with their argument is that Statement A is not 
true, as it would imply that we
would be able to write $J_{n,R}(g)=I_n(h)$ for some connection function $h$ with 
bounded support. This would mean that $g_{n,R}$ can be seen as a scaling of $h$, 
that is,
$$
1_{\{x \leq R\}} g(nx)=h(nx),
$$
but this leads to $h(x)=1_{\{x \leq nR\}}g(x)$, which clearly does not make any 
sense in general.

It seems then that Roy and Sarkar (2003) interchange truncation and scaling, but 
these two operations do not commute. This mixing up becomes already 
apparent when we look at 
their Lemma 5 which states (without proof) that 
\begin{eqnarray}
\lim_{n\to\infty} (\lambda_n \ell(K))^{-1} \IE L_{n,R}(g) &=& \plgr (1-\plgrr); 
\label{eq0d} \\
\lim_{n\to\infty} (\lambda_n \ell(K))^{-1} \Var L_{n,R}(g) &=& \plgr (1-\plgrr) 
+\lambda \int_{\IR^d} (1-g(|x|)) \nonumber \\
& & \Bigl[ \plgr^2 \plgrgr(x,0) - 2\plgr^2 \plgrr \plgrg(x,0) + \nonumber \\
& & \plg^2 \plgg(x,0) \Bigr] - \plgr^2 (1-\plgrr)^2\,dx; \label{eq0e}
\end{eqnarray}
where $\ell$ denotes Lebesgue measure on $\IR^d$ and
\[
p(\mu,h) = e^{-\mu \int_{\IR^d} h(|y|)\,dy} \quad \mbox{and} \quad
p_{\mu}^{h_1,h_2}(x_1,x_2) = e^{\mu \int_{\IR^d} h_1(|y-x_1|) h_2(|y-x_2|)\,dy}.
\]
However, the following proposition shows that (\ref{eq0d}) and ({\ref{eq0e}) are 
not 
correct; see the forthcoming Lemma \ref{corollary3} for a corresponding correct 
(and 
useful) statement. 
\begin{proposition} \sl
For\/\ $R>\sup\{|x_1-x_2|:x_1,x_2\in K\}$ we have
\begin{eqnarray}
\lim_{n\to\infty} (\lambda_n \ell(K))^{-1} \IE L_{n,R}(g) &=& 0 \label{eq0m}; \\
\lim_{n\to\infty} (\lambda_n \ell(K))^{-1} \Var L_{n,R}(g) &=& 0 \label{eq0n}.
\end{eqnarray}
\end{proposition}
{\bf Proof:} For $R>0$ and $n\in\IN$ define 
$k_{nR},k^{nR} 
: [0,\infty)\to [0,1]$ 
by
\[
k_{nR}(x)=1_{\{x\leq nR\}}g(x),\qquad k^{nR}(x)=1_{\{x>nR\}}g(x),\qquad x\in 
[0,\infty).
\]
We have as $n\to\infty$,
\[
\begin{array}{lclcl}
\plngnr &=& \plnndhnr &\to& \plg; \\
\plngnrr &=& \plnndhnnr &\to& 1.
\end{array}
\]
According to Roy and Sarkar (2003) Lemma 4 we have for $R>\sup\{|x_1-x_2|: 
x_1,x_2\in K\}$,
\begin{equation}
\IE L_R(g) = \lambda \ell(K) \plgr (1-\plgrr), \label{eq0b}
\end{equation}
and therefore,
\[
(\lambda_n\ell(K))^{-1} \IE L_{n,R}(g) = \plngnr (1-\plngnrr) \to 0,\qquad 
n\to\infty,
\]
which proves (\ref{eq0m}). 

To prove (\ref{eq0n}), we use Lemma 4 in Roy and Sarkar (2003) which says that 
for 
$R>\sup\{|x_1-x_2|: x_1,x_2\in K\}$, we have  
\begin{eqnarray}
\Var L_R(g) &=& \lambda \ell(K) \plgr (1-\plgrr) +\lambda^2 \int_K \int_K 
(1-g(|x_1-x_2|)) \nonumber \\
& & \Bigl[ \plgr^2 \plgrgr(x_1,x_2) - 2\plgr^2 \plgrr \plgrg(x_1,x_2) + 
\nonumber 
\\
& & \plg^2 \plgg(x_1,x_2) \Bigr] - \plgr^2 (1-\plgrr)^2 \,dx_2\,dx_1. 
\label{eq0c}
\end{eqnarray}
We use (\ref{eq0c}) with 
$\lambda=\lambda_n$ and $g=g_n$. Note that as $n\to\infty$
\[
\begin{array}{lclcl}
\plngnrgnr(x/n,0) &=& \plnndhnrhnr(x,0) &\to& \plgg(x,0); \\
\plngnrgn(x/n,0) &=& \plnndhnrg(x,0) &\to& \plgg(x,0); \\
\plngngn(x/n,0) &=& \plnndgg(x,0) &\to& \plgg(x,0).
\end{array}
\]
We have
\begin{eqnarray*}
\lefteqn{ \frac{\lambda_n}{\ell(K)} \int_K \int_K (1-g_n(|x_1-x_2|)) \Bigl[ 
\plngnr^2 \plngnrgnr(x_1,x_2) - } \\
\lefteqn{ 2\plngnr^2 \plngnrr \plngnrgn(x_1,x_2) + \plngn^2 \plngngn(x_1,x_2) 
\Bigr] - } \\
\lefteqn{ \plngnr^2 (1-\plngnrr)^2\,dx_2\,dx_1 = } \\
&=& \frac{\lambda_n}{n^d\ell(K)} \int_K \int_{n(K-x_1)} (1-g(|x_2|)) \Bigl[ 
\plnndhnr^2 \plnndhnrhnr(0,x_2) - \\
& & 2\plnndhnr^2 \plnndhnnr \plnndhnrg(0,x_2) + \plnndg^2 \plnndgg(0,x_2) \Bigr] 
- 
\\
& & \plnndhnr^2 (1-\plnndhnnr)^2\,dx_2\,dx_1.
\end{eqnarray*}
By Lemma \ref{lemma2} below with $x=-x_2$ we can apply the 
dominated 
convergence theorem. Combining the result with (\ref{eq0m}) yields (\ref{eq0n}). 
\QED

In what follows, we proceed along the way 
that we believe Roy and Sarkar (2003) had in mind.

For this, we introduce for $R>0$ and $n\in\IN$ connection functions $g_{R,n}, 
g^{R,n} : [0,\infty) \to [0,1]$ as follows:
\[
g_{R,n}(x)=1_{\{x\leq R/n\}}g(nx),\qquad g^{R,n}(x)=1_{\{x>R/n\}}g(nx).
\]
Note the difference between $g_{R,n}$ and $g_{n,R}$ and between $g^{R,n}$ and 
$g_n^R$.
Let $J_{R,n}(g)$ be the number of isolated vertices of $(X_n,\lambda_n, 
g_{R,n})$ 
in $K$ and let $L_{R,n}(g)=J_{R,n}(g)-I_n(g)$ be the number of isolated vertices 
of 
$(X_n,\lambda_n, g_{R,n})$ in $K$ that are not isolated in $(X_n,\lambda_n, 
g_n)$. 
Note that the notations $g_{R,n}$, $J_{R,n}(g)$ and $L_{R,n}(g)$ can formally 
not 
be used here, since $g_{n,R}$, $J_{n,R}(g)$ and $L_{n,R}(g)$ have already been 
defined. Nevertheless we shall adopt these notations, because henceforth we 
shall 
use the function $g_{n,R}$ and the random variables $J_{n,R}(g)$ and 
$L_{n,R}(g)$ 
no more. We now claim that the following is true (compare the incorrect 
Statement A 
above)

\medskip\noindent
{\bf Statement B}\, {\sl If (\ref{eq0}) is true when the connection function\/\ 
$g$ has bounded support, then it is the case that 
\begin{equation}\label{eq0p}
\frac{J_{R,n}(g)-\IE J_{R,n}(g)}{\sqrt{\Var J_{R,n}(g)}} \rightsquigarrow 
N(0,1),\qquad n\to\infty,
\end{equation}
for\/\ {\em any} connection function\/\ $g$.}

\medskip\noindent
To see this, observe that 
$$
J_{R,n}(g)=I_n(g_R),
$$
as can be seen by direct computation. Since $g_R$ has bounded support, Statement 
B 
follows. The moral of this is, that we should base the proof on $J_{R,n}(g)$ and 
$L_{R,n}(g)$ instead of $J_{n,R}(g)$ and $L_{n,R}(g)$. In the next section we 
show that 
the proof idea of Roy and Sarkar (2003) can still 
be carried out, although the computations involved are a little more complicated 
now. 

\section{Proof of Theorem \ref{mh}}

We start with a technical lemma, needed for applications of dominated 
convergence.

\begin{lemma}\label{lemma2} \sl There exists $N$ such that for $R>0$, $n\geq N$ 
and 
$x\in\IR^d$
\begin{eqnarray}
\lefteqn{ \Bigl| (1-g(|x|)) \Bigl[ \plnndgr^2 \plnndgrgr(x,0) - 2\plnndgr^2 
\plnndgrr \plnndgrg(x,0) + } \nonumber \\ 
& & \plnndg^2 \plnndgg(x,0) \Bigr] - \plnndgr^2 (1-\plnndgrr)^2 \Bigr| \leq C 
g(|x|/2),~~~~~ \label{eq7}
\end{eqnarray}
where $C$ is a constant not depending on $x$, $n$ or $R$.
\end{lemma}
{\bf Proof:} Since $\plnndgr \plnndgrr = \plnndg$, the expression between the 
absolute value signs in (\ref{eq7}) is equal to
\begin{eqnarray}
\lefteqn{ -g(|x|) \Bigl[ \plnndgr^2 \plnndgrgr(x,0) - 
2\plnndgr\plnndg\plnndgrg(x,0) + } \nonumber \\
& & \plnndg^2 \plnndgg(x,0) \Bigr] + \plnndgr^2 (\plnndgrgr(x,0)-1) - \nonumber 
\\
& & 2\plnndgr\plnndg (\plnndgrg(x,0)-1) + \plnndg^2 (\plnndgg(x,0)-1). 
\label{eq8}
\end{eqnarray}
Let $N$ be such that $\frac{3}{4}\lambda \leq \lambda_n/n^d \leq 
\frac{3}{2}\lambda$, $n\geq N$. Then since
\[
\int_{\IR^d} g_R(|y|)+g(|y|)\,dy \geq 2 \int_{\IR^d} g_R(|y|)\,dy \geq 2 
\int_{\IR^d} g_R(|y|)g(|y+x|)\,dy,
\]
we have for $n\geq N$
\begin{equation}\label{eq9}
\plnndgr\plnndg\plnndgrg(x,0) \leq e^{-\frac{3}{4}\lambda \int_{\IR^d} 
g_R(|y|)+g(|y|)\,dy + \frac{3}{2}\lambda \int_{\IR^d} g_R(|y-x|)g(|y|)\,dy} \leq 
1.
\end{equation}
Also,
\[
\plnndgr^2\plnndgrgr(x,0) \leq 1,\qquad \plnndg^2\plnndgg(x,0) \leq 1,
\]
which follows from (\ref{eq9}) by taking $g=g_R$ or letting $R\to\infty$ 
respectively. Hence for $n\geq N$ the absolute value of (\ref{eq8}) is bounded 
by
\begin{equation}\label{eq10}
4g(|x|) + 4(p_{2\lambda}^{g,g}(x,0)-1).
\end{equation}

To give an upper bound for the second term in this expression, note that for 
$y\in\IR^d$ either $|y|\geq|x|/2$ or $|y-x|\geq |x|/2$, so
\begin{eqnarray*}
\int_{\IR^d} g(|y-x|)g(|y|)\,dy &\leq& \int_{|y|<|x|/2} g(|y-x|)g(|y|)\,dy + 
\int_{|y|\geq |x|/2} g(|y-x|)g(|y|)\,dy \\
&\leq& g(|x|/2) \int_{|y|<|x|/2} g(|y|)\,dy + g(|x|/2) \int_{|y|\geq |x|/2} 
g(|y-x|)\,dy \\
&\leq& 2 g(|x|/2) \int_{\IR^d} g(|y|)\,dy.
\end{eqnarray*}
Choose $M$ such that $4\lambda g(M/2) \int_{\IR^d} g(|y|)\,dy \leq 1$. Then 
since 
$e^t \leq 1+et$, $t\leq 1$, we have for $|x|\geq M$
\[
e^{4\lambda g(|x|/2) \int_{\IR^d} g(|y|)\,dy} \leq 1+4e\lambda 
g(|x|/2)\int_{\IR^d} 
g(|y|)\,dy.
\]
For $|x|<M$ we have
\[
e^{4\lambda g(|x|/2) \int_{\IR^d} g(|y|)\,dy} \leq e^{4\lambda \int_{\IR^d} 
g(|y|)\,dy} \leq 1 + g(|x|/2) g(M/2)^{-1} [e^{4\lambda \int_{\IR^d} g(|y|)\,dy} 
-1].
\]
Combining the above inequalities yields
\begin{equation}\label{eq10a}
p_{2\lambda}^{g,g}(x,0) -1 \leq C g(|x|/2),
\end{equation}
where $C$ is a constant not depending on $x$, $n$ or $R$. We conclude that 
(\ref{eq10}) is bounded by $4(1+C)g(|x|/2)$. \QEDnoskip

\begin{lemma}\label{lemma1} \sl \ 
\begin{eqnarray}
\IE L_{R,n}(g) &=& \lambda_n \ell(K) \plngrn (1-\plngrrn) \label{eq1} \\
\Var L_{R,n}(g) &=& \lambda_n \ell(K) \plngrn (1-\plngrrn) +\lambda_n^2 \int_K 
\int_K (1-g_n(|x_1-x_2|)) \nonumber \\
& & \Bigl[ \plngrn^2 \plngrngrn(x_1,x_2) - 2\plngrn^2 \plngrrn 
\plngrngn(x_1,x_2) + 
\nonumber \\
& & \plngn^2 \plngngn(x_1,x_2) \Bigr] - \plngrn^2 (1-\plngrrn)^2 \,dx_2\,dx_1 + 
\nonumber \\
& & \lambda_n^2 \plngrn^2 \int_K \int_K g^{R,n}(|x_1-x_2|) 
\plngrngrn(x_1,x_2)\,dx_2\,dx_1 \label{eq2}
\end{eqnarray}
\end{lemma}
{\bf Proof:} The first statement (\ref{eq1}) is proved as in Roy and Sarkar 
(2003) 
Lemma 4.

For a Borel subset $B$ of $\IR^d$ let $X_n(B)$ be the number of 
points 
in $X_n\cap B$. For $t>0$ denote $K^t=K+\{x\in\IR^d : |x|< t\}$. In the model 
$(X_n,\lambda_n, g_n)$ let $L_{R,n,t}(g)$ be the number of points $\xi$ in 
$X_n\cap 
K$ such that $\xi$ is not connected to any point in $X_n\cap K^t$ at a distance 
$R/n$ or less from $\xi$ but $\xi$ is connected to some point in $X_n\cap K^t$ 
at a 
distance greater than $R/n$ from $\xi$. Since $L_{R,n,t}(g) \to L_{R,n}(g)$, 
$t\to\infty$, and $L_{R,n,t}(g)\leq X_n(K)$, $t>0$, and $\IE X_n(K)^2 < \infty$, 
the dominated convergence theorem gives
\[
\IE L_{R,n,t}(g) \to \IE L_{R,n}(g), \qquad \Var L_{R,n,t}(g) \to \Var 
L_{R,n}(g), 
\qquad t\to\infty.
\]
In order to compute the moments of $L_{R,n,t}(g)$, note that
\[
L_{R,n,t}(g) \sim \sum_{i=1}^{X_n(K^t)} 1_{F_i},
\]
where $\sim$ denotes equality in distribution, $\xi_i$, $i\in\IN$ are 
independent 
random variables, independent of $X_n(K^t)$, uniformly distributed on $K^t$ and 
connected to each other according to $g_n$, and $F_i = \{\xi_i\in K$; $\xi_i$ is 
not connected to any $\xi_j$, $j\leq X_n(K^t)$, at a distance $R/n$ or less from 
$\xi_i$; $\xi_i$ is connected to some $\xi_j$, $j\leq X_n(K^t)$, at a distance 
greater than $R/n$ from $\xi_i \}$.

Since
\[
L_{R,n,t}(g)^2 \sim \sum_{i=1}^{X_n(K^t)} 1_{F_i} + \sum_{i=1}^{X_n(K^t)} 
\sum_{{\scriptstyle j=1}\atop{\scriptstyle j\not= i}}^{X_n(K^t)} 1_{F_i} 
1_{F_j},
\]
the variance of $L_{R,n,t}(g)$ can be written as
\begin{eqnarray}
\lefteqn{ \Var L_{R,n,t}(g) = } \label{eq3} \\
&\!\! =& \!\! \IE L_{R,n,t}(g) + \sum_{m=2}^{\infty} m(m-1) \IP(F_1\cap F_2\,|\, 
X_n(K^t)=m) \IP(X_n(K^t)=m) - (\IE L_{R,n,t}(g))^2. \nonumber
\end{eqnarray}

We have
\begin{eqnarray*}
\lefteqn{ \IP(F_1\cap F_2\cap \{\xi_1\mbox{ is connected to }\xi_2\}\,|\, 
X_n(K^t)=m) = } \\
&=& \frac{1}{\ell(K^t)^m} \int_K \int_K g^{R,n}(|x_1-x_2|) \int_{K^t} \ldots 
\int_{K^t} \\
& & \prod_{i=3}^m (1-g_{R,n}(|x_i-x_1|))(1-g_{R,n}(|x_i-x_2|) \,dx_m\ldots 
dx_3\,dx_2\,dx_1 \\
&=& \frac{1}{\ell(K^t)^m} \int_K \int_K g^{R,n}(|x_1-x_2|) \\
& & \biggl[ \int_{K^t} (1-g_{R,n}(|y-x_1|))(1-g_{R,n}(|y-x_2|))\,dy 
\biggr]^{m-2} 
dx_2\,dx_1,
\end{eqnarray*}
whence
\begin{eqnarray}
\lefteqn{ \sum_{m=2}^{\infty} m(m-1) \IP(F_1\cap F_2\cap \{\xi_1\mbox{ is 
connected 
to }\xi_2\}\,|\, X_n(K^t)=m) \IP(X_n(K^t)=m) =} \nonumber \\
&=& \lambda_n^2 \int_K \int_K g^{R,n}(|x_1-x_2|) \sum_{m=0}^{\infty} 
\frac{e^{-\lambda_n \ell(K^t)} \lambda_n^m}{m!} \nonumber \\
& & \biggl[ \int_{K^t} (1-g_{R,n}(|y-x_1|))(1-g_{R,n}(|y-x_2|))\,dy \biggr]^m 
dx_2\,dx_1 \nonumber \\
&=& \lambda_n^2 \int_K \int_K g^{R,n}(|x_1-x_2|) e^{\lambda_n \int_{K^t} 
-g_{R,n}(|y-x_1|) - g_{R,n}(|y-x_2|) + g_{R,n}(|y-x_1|)g_{R,n}(|y-x_2|)\,dy} 
dx_2\,dx_1 \nonumber \\
&\to& \lambda_n^2 \plngrn^2 \int_K \int_K g^{R,n}(|x_1-x_2|) 
\plngrngrn(x_1,x_2)\,dx_2\,dx_1,\qquad t\to\infty, \label{eq4}
\end{eqnarray}
where we use the dominated convergence theorem.

Furthermore,
\begin{eqnarray}
\lefteqn{ \IP(F_1\cap F_2\cap \{\xi_1\mbox{ is not connected to }\xi_2\}\,|\, 
X_n(K^t)=m) = } \nonumber \\
&=& \frac{1}{\ell(K^t)^m} \int_K \int_K (1-g_n(|x_1-x_2|)) \int_{K^t} \ldots 
\int_{K^t} \nonumber \\
& & \biggl[ 1-\prod_{i=3}^m (1-g^{R,n}(|x_i-x_1|)) \biggr] \prod_{i=3}^m 
(1-g_{R,n}(|x_i-x_1|)) \nonumber \\
& & \biggl[ 1-\prod_{i=3}^m (1-g^{R,n}(|x_i-x_2|)) \biggr] \prod_{i=3}^m 
(1-g_{R,n}(|x_i-x_2|)) \,dx_m\ldots dx_3\,dx_2\,dx_1. \label{eq5}
\end{eqnarray}
Exactly as in Roy and Sarkar (2003) Lemma 4, one can now show that
\begin{eqnarray}
\lefteqn{ \sum_{m=2}^{\infty} m(m-1) \IP(F_1\cap F_2\cap \{\xi_1\mbox{ is not 
connected to }\xi_2\}\,|\, X_n(K^t)=m) \IP(X_n(K^t)=m) =} \nonumber \\
&\to& \lambda_n^2 \int_K \int_K (1-g_n(|x_1-x_2|)) \Bigl[ \plngrn^2 
\plngrngrn(x_1,x_2) - \nonumber \\
& & 2\plngrn^2 \plngrrn \plngrngn(x_1,x_2) + \plngn^2 \plngngn(x_1,x_2) \Bigr] 
dx_2\,dx_1, \label{eq6}
\end{eqnarray}
as $t\to\infty$, where we use the dominated convergence theorem.

Combining (\ref{eq3}), (\ref{eq1}), (\ref{eq4}) and (\ref{eq6}) yields 
(\ref{eq2}). 
\QED

The following lemma replaces the incorrect Lemma 5 (equation (\ref{eq0d}) and 
(\ref{eq0e}) in our current paper) of Roy and Sarkar (2003).

\begin{lemma}\label{corollary3} \sl \ 
\begin{eqnarray}
\lim_{n\to\infty} (\lambda_n \ell(K))^{-1} \IE L_{R,n}(g) &=& \plgr (1-\plgrr) 
\label{eq11} \\
\lim_{n\to\infty} (\lambda_n \ell(K))^{-1} \Var L_{R,n}(g) &=& \plgr (1-\plgrr) 
+\lambda \int_{\IR^d} (1-g(|x|)) \nonumber \\
& & \Bigl[ \plgr^2 \plgrgr(x,0) - 2\plgr^2 \plgrr \plgrg(x,0) + \nonumber \\
& & \plg^2 \plgg(x,0) \Bigr] - \plgr^2 (1-\plgrr)^2\,dx + \nonumber \\
& & \lambda \plgr^2 \int_{\IR^d} g^R(|x|) \plgrgr(x,0)\,dx \label{eq12}
\end{eqnarray}
\end{lemma}
{\bf Proof:} Assertion (\ref{eq11}) follows from (\ref{eq1}) by direct 
computation. 
We shall deduce (\ref{eq12}) from (\ref{eq2}). By the dominated convergence 
theorem
\begin{eqnarray}
\lefteqn{ \frac{\lambda_n}{\ell(K)} \plngrn^2 \int_K \int_K g^{R,n}(|x_1-x_2|) 
\plngrngrn(x_1,x_2)\,dx_2\,dx_1 = } \nonumber \\
&=& \frac{\lambda_n}{n^d \ell(K)} \plngrn^2 \int_K \int_{n(K-x_1)} g^R(|x_2|) 
\plnndgrgr(0,x_2)\,dx_2\,dx_1 \nonumber \\
&\to& \lambda \plgr^2 \int_{\IR^d} g^R(|x|) \plgrgr(x,0)\,dx,\qquad n\to\infty. 
\label{eq13}
\end{eqnarray}
Furthermore,
\begin{eqnarray*}
\lefteqn{ \frac{\lambda_n}{\ell(K)} \int_K \int_K (1-g_n(|x_1-x_2|)) \Bigl[ 
\plngrn^2 \plngrngrn(x_1,x_2) - } \\
\lefteqn{ 2\plngrn^2 \plngrrn \plngrngn(x_1,x_2) + \plngn^2 \plngngn(x_1,x_2) 
\Bigr] - } \\
\lefteqn{ \plngrn^2 (1-\plngrrn)^2 \,dx_2\,dx_1 = } \\
&=& \frac{\lambda_n}{n^d\ell(K)} \int_K \int_{n(K-x_1)} (1-g(|x_2|)) \Bigl[ 
\plnndgr^2 \plnndgrgr(0,x_2) - \\
& & 2\plnndgr^2\plnndgrr\plnndgrg(0,x_2) + \plnndg^2\plnndgg(0,x_2) \Bigr] - \\
& & \plnndgr^2 (1-\plnndgrr)^2\,dx_2\,dx_1.
\end{eqnarray*}
By Lemma \ref{lemma2} with $x=-x_2$, we can apply the dominated convergence 
theorem. Combining the result with (\ref{eq11}) and (\ref{eq13}) yields 
(\ref{eq12}). \QEDnoskip

\begin{corollary}\label{corollary4} \sl \ 
\begin{eqnarray}
\lim_{R\to\infty} \lim_{n\to\infty} (\lambda_n \ell(K))^{-1} \IE L_{R,n}(g) &=& 
0 
\label{eq14} \\
\lim_{R\to\infty} \lim_{n\to\infty} (\lambda_n \ell(K))^{-1} \Var L_{R,n}(g) &=& 
0 
\label{eq15}
\end{eqnarray}
\end{corollary}
{\bf Proof:} The dominated convergence theorem gives
\[
\plgr\to\plg,\quad \plgrr\to 1,\quad \plgrgr(x,0)\to \plgg(x,0),\quad 
\plgrg(x,0)\to\plgg(x,0),
\]
as $R\to\infty$. Now (\ref{eq14}) follows from (\ref{eq11}). Another application 
of 
the dominated convergence theorem yields
\[
\int_{\IR^d} g^R(|x|) \plgrgr(x,0)\,dx \to 0,\qquad R\to\infty.
\]
Finally, the integrand in the first integral on the right hand side of 
(\ref{eq12}) 
tends to $0$ as $R\to\infty$. By Lemma \ref{lemma2} with $\lambda_n=\lambda 
n^d$, 
we can apply the dominated convergence theorem to conclude (\ref{eq15}). \QED

Finally, we can prove the main result:

\begin{theorem}\label{theorem5} \sl If for $R>0$
\begin{equation}\label{eq16}
\frac{J_{R,n}(g)-\IE J_{R,n}(g)}{\sqrt{\Var J_{R,n}(g)}} \rightsquigarrow 
N(0,1),\qquad n\to\infty,
\end{equation}
then (\ref{eq0}) holds.
\end{theorem}
{\bf Proof:} Roy and Sarkar (2003) Lemma 3 shows that
\begin{equation}\label{eq0h}
\lim_{n\to\infty} (\lambda_n \ell(K))^{-1} \Var I_n(g) = \plg + \lambda \plg^2 
\int_{\IR^d} (1-g(|x|)) \plgg(x,0) -1 \,dx.
\end{equation}
It follows from (\ref{eq0h}), Corollary \ref{corollary4} and 
Chebyshev's inequality that
\[
\lim_{R\to\infty} \limsup_{n\to\infty} \IP \biggl( \biggl| \frac{L_{R,n}(g)-\IE 
L_{R,n}(g)}{\sqrt{\Var I_n(g)}} \biggr| \geq \varepsilon \biggr) \leq 
\lim_{R\to\infty} \lim_{n\to\infty} \frac{\Var L_{R,n}(g)}{\varepsilon^2 \Var 
I_n(g)} = 0,\qquad \varepsilon >0.
\]
Moreover, applying (\ref{eq0h}) also with $g$ replaced by $g_R$ gives 
$\lim_{n\to\infty} \Var 
J_{R,n}(g) / \Var I_n(g) = \delta_R$, where $\delta_R$ is a constant. 
(This was incorrectly claimed in Roy and Sarkar (2003) with $L_{n,R}(g)$ instead 
of 
$L_{R,n}(g)$.) Because
\[
(1-g_R(|x|)) \plgrgr(x,0) -1 \geq (1-g_R(|x|))\cdot 1 -1 \geq -g(|x|)
\]
and by (\ref{eq10a})
\[
(1-g_R(|x|)) \plgrgr(x,0) -1 \leq 1\cdot \plgg(x,0) -1 \leq Cg(|x|/2),
\]
where $C$ is a constant not depending on $x$ or $R$, we have by the dominated 
convergence theorem $\lim_{R\to\infty} \delta_R =1$. Now if (\ref{eq16}) holds, 
then for $x\in\IR$
\begin{eqnarray*}
\lefteqn{ \limsup_{n\to\infty} \IP\biggl( \frac{I_n(g)-\IE I_n(g)}{\sqrt{\Var 
I_n(g)}} \leq x \biggr) \leq } \\
&\leq& \lim_{\varepsilon\downarrow 0} \lim_{R\to\infty} \limsup_{n\to\infty} 
\IP\biggl( \frac{J_{R,n}(g)-\IE J_{R,n}(g)}{\sqrt{\Var I_n(g)}} \leq x 
+\varepsilon 
\biggr) + \IP \biggl( \biggl| \frac{L_{R,n}(g)-\IE L_{R,n}(g)}{\sqrt{\Var 
I_n(g)}} 
\biggr| \geq \varepsilon \biggr) \\
&=& \Phi(x).
\end{eqnarray*}
A similar argument yields
\[
\liminf_{n\to\infty} \IP\biggl( \frac{I_n(g)-\IE I_n(g)}{\sqrt{\Var I_n(g)}} 
\leq x 
\biggr) \geq \Phi(x),
\]
which completes the proof of the theorem. \QEDnoskip

\section{Extension to larger componenents}

In this section, we discuss larger components. A central limit theorem for
larger components needs another approach, even when the connection function has
bounded support. The reason for this is that the exact moment computations of 
the preceding sections no longer seem possible. At this point, we can only prove
a central limit theorem when the connection function $g$ has bounded support. 
For this, we use a result of Bolthausen (1982), from which it follows that in 
order to prove a central limit theorem, certain mixing conditions suffice. For 
convenience, the central limit theorem in this section
is stated a little different from the earlier ones, in the sense that we do not
scale the connection function and the density, but instead take larger and
larger subsets of the space. This is equivalent to the case where 
$\lambda_n=\lambda n^d$ in the original setup.

For a subset $\Lambda$ of $\IZ^d$, let the inner boundary of $\Lambda$ be 
denoted by $\partial \Lambda$, and its cardinality by $|\Lambda|$. Let the 
random variable $I^r(\Lambda)=I^r(\Lambda,g)$ be defined as $1/r$ times the 
number of vertices 
of $(X,\lambda,g)$ in $\Lambda+(0,1]^d$ that are contained in a component of 
size 
$r$. For $z\in\IZ^d$ write $I^r(z)=I^r(\{z\})$. We shall prove the following 
central limit theorem.

\begin{theorem}
\label{ncrv}
 \sl
Consider a random connection model with connection function\/\ $g$ of bounded 
support. Then for any increasing sequence\/\ $(\Lambda_n)_{n\in\IN}$ of finite 
subsets of\/\ $\IZ^d$ with\/\ $\bigcup_{n\in\IN} \Lambda_n = \IZ^d$ and\/\ 
$|\partial\Lambda_n|/|\Lambda_n|\to 0$, $n\to\infty$, we have
\begin{equation}\label{eq28fa}
\frac{I^r(\Lambda_n)-\IE I^r(\Lambda_n)}{\sqrt{\Var I^r(\Lambda_n)}} 
\rightsquigarrow N(0,1),\qquad n\to\infty.
\end{equation}
\end{theorem}

In order to prove this result, we use the main theorem in Bolthausen (1982). The 
conditions of his theorem involve three mixing conditions which are trivially
satisfied when $g$ has bounded support, and which we do not repeat here. Under 
these three mixing conditions, Bolthausen (1982) shows that if in addition 
\begin{equation}\label{eq28e}
\sum_{z\in\IZ^d} \Cov(I^r(0),I^r(z)) > 0,
\end{equation}
then it is the case that
\begin{equation}\label{eq28f}
\frac{I^r(\Lambda_n)-\IE I^r(\Lambda_n)}{\sqrt{|\Lambda_n| \sum_{z\in\IZ^d} 
\Cov(I^r(0),I^r(z))}} \rightsquigarrow N(0,1),\qquad n\to\infty.
\end{equation}
Because of the following elementary lemma, which we give without proof, 
(\ref{eq28e}) and (\ref{eq28f}) imply our Theorem \ref{ncrv}.

\begin{lemma}\label{lemma7} \sl
Let\/\ $(Y_z)_{z\in\IZ^d}$ be a stationary random field with\/\ $\IE 
Y_0^2<\infty$. Let\/\ $(\Lambda_n)_{n\in\IN}$ be a sequence of finite non-empty 
subsets of\/\ $\IZ^d$ with\/\ $|\partial\Lambda_n|/|\Lambda_n|\to 0$, 
$n\to\infty$. If
\begin{equation}\label{sumcov}
\sum_{z\in\IZ^d} |\Cov(Y_0,Y_z)| < \infty,
\end{equation}
then
\[
\frac{1}{|\Lambda_n|} \Var \sum_{z\in\Lambda_n} Y_z \to \sum_{z\in\IZ^d} 
\Cov(Y_0,Y_z),\qquad n\to\infty.
\]
\end{lemma}

\noindent
Note that (\ref{sumcov}) is satisfied because $g$ has bounded support. It 
remains to prove (\ref{eq28e}). We give the proof in the 
two-dimensional case, but the method clearly generalizes to other dimensions.

With a slight abuse of notation, for a Borel subset $B$ of $\IR^2$ let $I^r(B)$ 
henceforth be defined as $1/r$ times the number of vertices of $(X,\lambda,g)$ 
in $B$ that are contained in a component of size $r$.
According to  Lemma \ref{lemma7}, it 
suffices to show that there exists 
$M\in\IN$ and $\gamma >0$ such that for all $n$,
\begin{equation}\label{eq3001}
{\Var I^r((0,nM]^2)} \geq \gamma n^2.
\end{equation}
We estimate the variance in (\ref{eq3001}) with the following general abstract 
trick, which we learned from J.\ v.d.\ Berg (personal communication).

\begin{lemma}\label{lemma301} \sl
Let\/\ $Y$ be a random variable with finite second moment, defined on a 
probability space\/\ $(\Omega,\AAA,\IP)$. Let\/\ $n\in\IN$ and let\/\ 
$\FFF_0\subseteq\FFF_1\subseteq\cdots\subseteq\FFF_n$ be sub-$\sigma$-algebras 
of\/\ $\AAA$ with\/\ $\IE(Y|\FFF_0)=\IE Y$ and\/\ $\IE(Y|\FFF_n)=Y$ a.s. Then we 
have
\[
\Var Y = \sum_{i=1}^n \IE\left[\IE(Y|\FFF_i)-\IE(Y|\FFF_{i-1})\right]^2.
\]
\end{lemma}
{\bf Proof:} For $1 \leq i \leq n$, denote 
$\Delta_i=\IE(Y|\FFF_i)-\IE(Y|\FFF_{i-1})$. We write the variance of $Y$ with a 
telescoping sum as $\Var Y=\IE(\sum_{i=1}^n 
\Delta_i)^2$. For $1\leq i<j\leq n$ we have $\IE\Delta_i\Delta_j=\IE\IE(\Delta_i 
Y|\FFF_j)-\IE\IE(\Delta_i Y|\FFF_{j-1})=0$. Hence $\Var Y=\sum_{i=1}^n 
\IE\Delta_i^2$, as required. \QED

Let $R$ be such that $g(x)=0$, $x\geq R$. Define $\mu=\IE I^r((0,1]^2)>0$. 
Choose an integer $M>3\lambda R/\mu$. We shall show that (\ref{eq3001}) holds 
for this $M$, and this is sufficient to prove Theorem \ref{ncrv}.

Partition the first quadrant of $\IR^2$ into cubes of side length $M$, and 
denote these cubes by $B_k$, $k\in\IN$, where the indices run 
as indicated in Figure 1. For $n\in\IN$ let $K_n$ be the set of indices 
$k\in\{1,\ldots,(n-1)^2\}$ that are shaded in Figure 1.

For $k\in \bigcup_{n \in \IN}K_n$, we define the following sets:
\begin{eqnarray*}
\verschovenbox_k &=& (rR,rR)+B_k; \\
\grotebox_k &=& B_k + (-rR,rR]^2; \\
\gebiedL_k &=& \textstyle \grotebox_k \cap \bigcup_{i=1}^{k-1} B_i; \\
\gebiedgeroteerdeL_k &=& \textstyle \grotebox_k \setminus \bigcup_{i=1}^{k-1} 
B_i;
\end{eqnarray*}
see Figure 2 and Figure 3.

\begin{center}
\scalebox{0.6}{\includegraphics{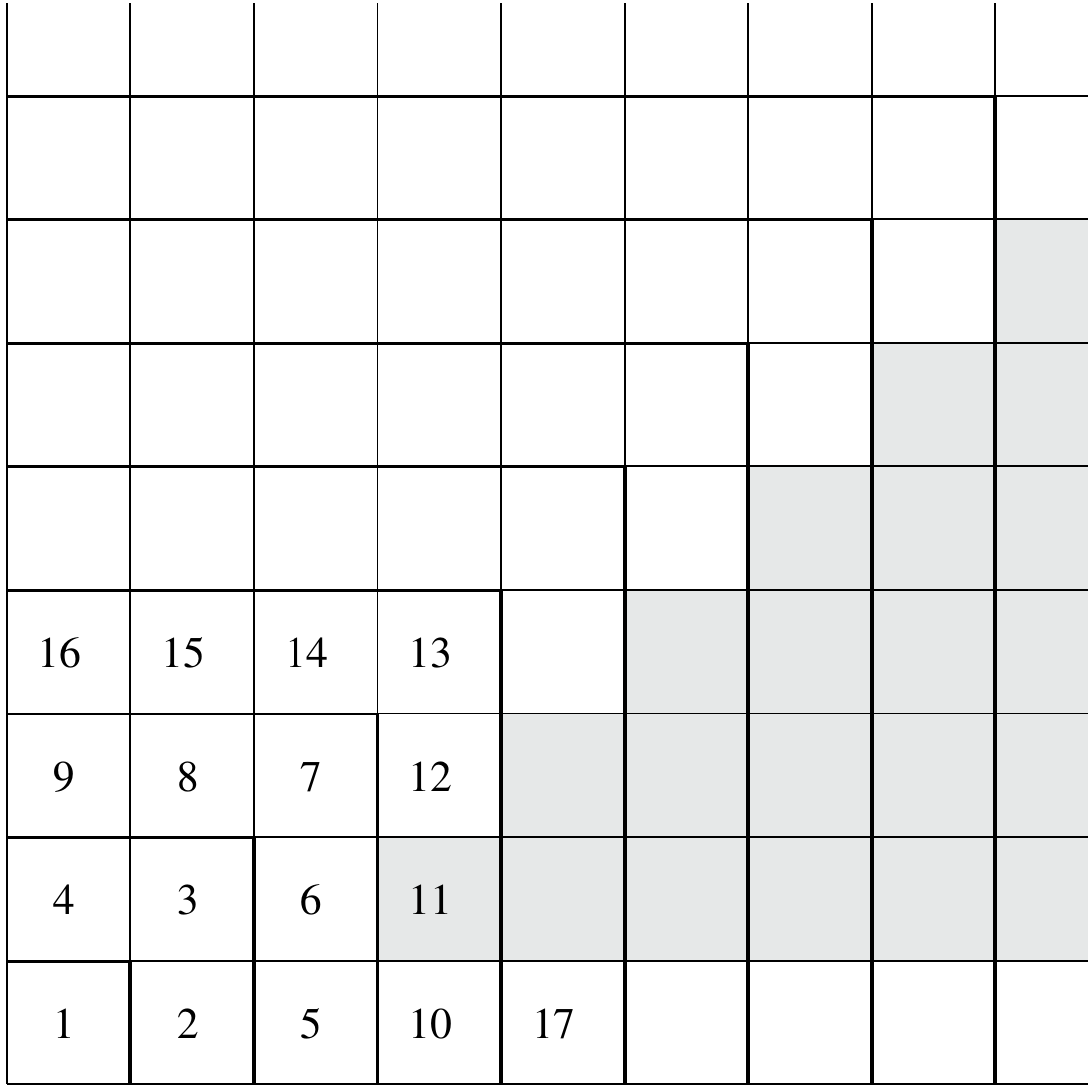}}
\end{center}
\begin{center}
Figure 1: The enumeration of cubes in the first quadrant.
\end{center}

\begin{center}
\scalebox{0.6}{\includegraphics{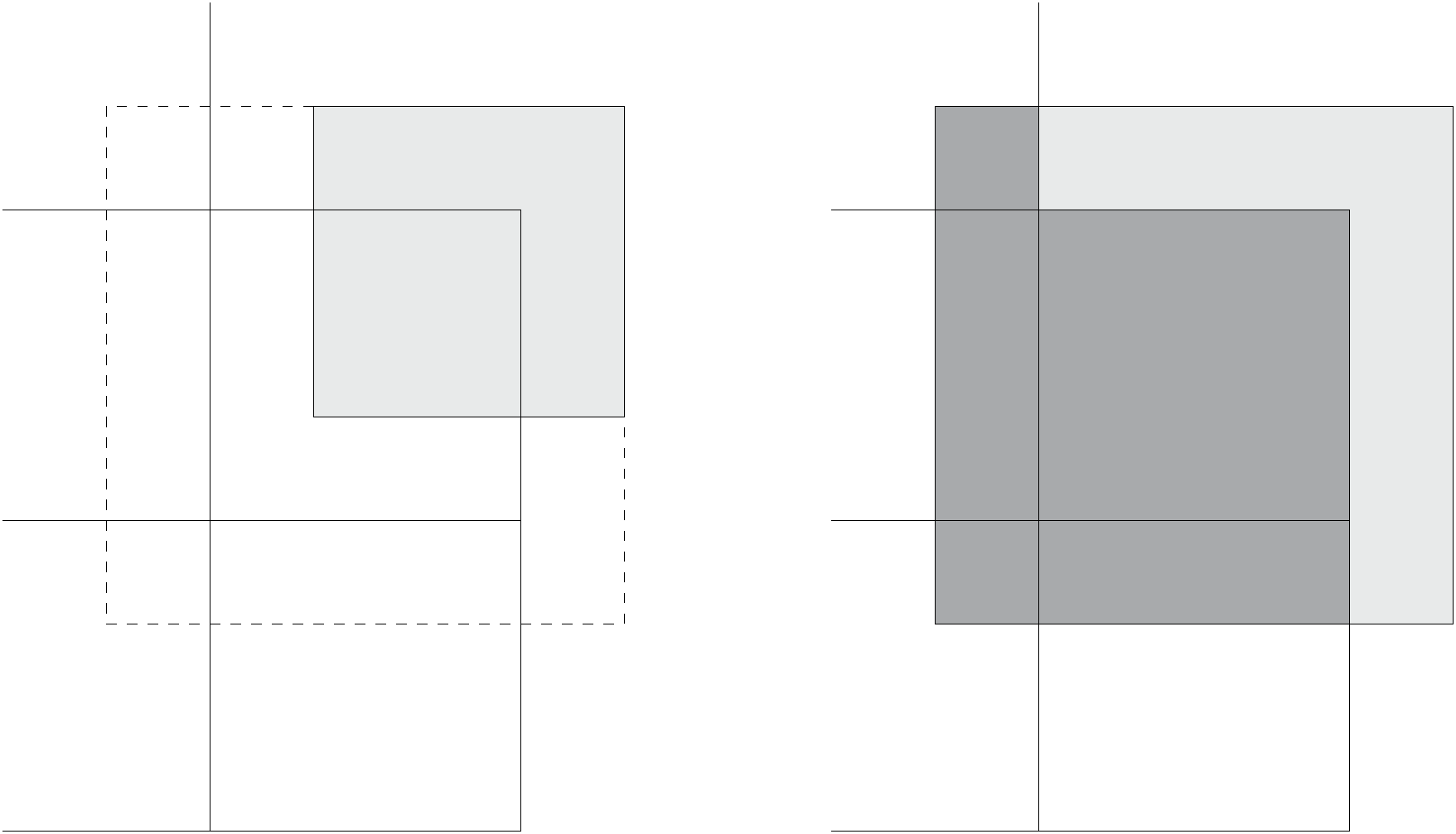}}
\end{center}
\begin{center}
Figure 2\hspace*{4.2cm} Figure 3
\end{center}
\begin{center} 
The shaded region on the left is $\verschovenbox_k$. The dark shaded region on 
the right is $\gebiedL_k$ and the light shaded region on the right is 
$\gebiedgeroteerdeL_k$.
\end{center}

For $k\in\IN$, let $\FFF_k$ be the $\sigma$-algebra generated by the points
of $X$ in $\bigcup_{i=1}^k B_i$. We shall first show that for $n\in\IN$ and 
$k\in K_n$ the difference 
$\IE(I^r((0,nM]^2)|\FFF_{k-1})-\IE(I^r((0,nM]^2)|\FFF_k)$ is bounded below by a 
positive uniform constant, with positive probability which is also uniform in 
$k$ and $n$.

On the one hand, we have
\begin{eqnarray}
\IE(I^r((0,nM]^2)|\FFF_{k-1}) &\geq& \IE(I^r(\verschovenbox_k)|\FFF_{k-1}) + 
\IE(I^r((0,nM]^2\setminus\grotebox_k)|\FFF_{k-1}) \nonumber \\
&=& \mu M^2 + \IE(I^r((0,nM]^2\setminus\grotebox_k)|\FFF_k), \label{eq307}
\end{eqnarray}
since $I^r(\verschovenbox_k)$ is independent of $\FFF_{k-1}$ and since the 
$\sigma$-algebra generated by $I^r((0,nM]^2\setminus\grotebox_k)$ and the points 
of $X$ in $\bigcup_{i=1}^{k-1}B_i$, is independent of the points of $X$ in 
$B_k$. 

On the other hand, we also have
\begin{eqnarray}
\lefteqn{ \IE(I^r((0,nM]^2)|\FFF_k) \leq } \nonumber \\
&\leq& (1/r)\IE(X(\gebiedL_k)|\FFF_k) + (1/r)\IE(X(\gebiedgeroteerdeL_k)|\FFF_k) 
+ \IE(I^r((0,nM]^2\setminus\grotebox_k)|\FFF_k) \nonumber \\
&=& 0 + 2\lambda R(M+rR) + \IE(I^r((0,nM]^2\setminus\grotebox_k)|\FFF_k), 
\label{eq308}
\end{eqnarray}
with probability at least $e^{-\lambda(M+2rR)^2}$, since $X(\gebiedL_k)$ is 
$\FFF_k$-measurable and $X(\gebiedgeroteerdeL_k)$ is independent of $\FFF_k$.

Combining (\ref{eq307}) and (\ref{eq308}) yields for $n\in\IN$ and $k\in K_n$,
\[
\IP\left(\IE(I^r((0,nM]^2)|\FFF_{k-1})-\IE(I^r((0,nM]^2)|\FFF_k)\geq \mu 
M^2-2\lambda R(M+rR)\right) \geq e^{-\lambda(M+2rR)^2}.
\]
Now observe that the box $(0,nM]^2$ contains at least $\alpha n^2$ boxes indexed
by an element of $K_n$, for some $\alpha >0$. Hence, 
since $\mu M^2-2\lambda R(M+rR)>0$, we have by Lemma \ref{lemma301}
\begin{eqnarray*}
\Var I^r((0,nM]^2) &\geq& \sum_{k\in K_n} \IE[ \IE(I^r((0,nM]^2)|\FFF_k) - 
\IE(I^r((0,nM]^2)|\FFF_{k-1}) ]^2 \\
&\geq& \alpha n^2 (\mu M^2-2\lambda R(M+rR))^2 e^{-\lambda(M+2rR)^2}, 
\end{eqnarray*}
proving the result.

\section*{References}

\medskip\noindent
Bolthausen, E. (1982), {\em On the central limit theorem for stationary
mixing random fields\/}, The Annals of Probability {\bf 10}, 1047--1050.

\medskip\noindent
Roy, R. and Sarkar, A. (2003), {\em High density asymptotics of the Poisson 
random 
connection model\/}, Physica A {\bf 318}, 230--242.
\end{document}